\documentclass[12pt,draft]{amsart}
\usepackage[all]{xy}
\usepackage{amsfonts}
\usepackage{amssymb}
\newtheorem{teo}{Theorem}[section]
\newtheorem{lem}[teo]{Lemma}
\newtheorem{prop}[teo]{Proposition}

\newtheorem{dfn}[teo]{Definition}
\newtheorem{rk}[teo]{Remark}
\newtheorem{ex}[teo]{Example}

\def\<{\langle}
\def\>{\rangle}

\def\e{\varepsilon}

\def\r{\rho}

\def\f{{\varphi}}

\def\C{{\mathbb C}}

\def\Z{{\mathbb Z}}

\def\End{\mathop{\rm End}\nolimits}

\def\Aut{\mathop{\rm Aut}\nolimits}

\def\Id{\mathop{\rm Id}\nolimits}
\def\Tr{\mathop{\rm tr}\nolimits}

\def\1{\mathbf 1}

\newcommand{\ov}[1]{\overline{#1}}

\newcommand{\wh}[1]{\widehat{#1}}

\def\Fix{\mathop{\rm Fix}\nolimits}
\def\Coin{\mathop{\rm Coin}\nolimits}

\begin{document}

\title[Bitwisted Burnside-Frobenius and Dehn]
{Bitwisted  Burnside-Frobenius theorem  and  Dehn conjugacy problem}
\author{Alexander Fel'shtyn}
\address{Instytut Matematyki, Uniwersytet Szczecinski,
ul. Wielkopolska 15, 70-451 Szczecin, Poland
 and  Boise State University, 1910
University Drive, Boise, Idaho, 83725-155, USA,}
\email{:felshtyn@diamond.boisestate.edu}


\begin{abstract}
It is proved for  Abelian groups  that the Reidemeister
coincidence number of two endomorphisms $\phi$ and $\psi$
is equal to the number of  coincidence  points of $\wh\phi$
and $\wh\psi$  on the unitary dual, if the  Reidemeister number  is finite.
An affirmative
answer to the bitwisted Dehn conjugacy problem
for almost polycyclic  groups  is obtained.
Finally we explain why the Reidemeister numbers are always infinite
for injective endomorphisms of  Baumslag-Solitar groups.
\end{abstract}

\maketitle
\tableofcontents

\section{Introduction}

\begin{dfn}
Let $G$ be a countable discrete group  and $\phi, \psi : G\rightarrow G$
two endomorphisms. Two elements $x,x'\in G$ are said to be
$(\phi, \psi)-conjugate$ iff there exists $\gamma \in G$ with
$$
x'=\psi(\gamma)  x   \phi(\gamma)^{-1}.
$$
The number of $(\phi,\psi$)-conjugacy classes is called the Reidemeister coincidence
number of   endomorphisms $\phi$ and $\psi$, denoted by $R(\phi,\psi)$.  If $\psi$ is the identity map then the $(\phi, id) $-conjugacy classes are the $\phi$- conjugacy classes in the group $G$ and $R(\phi,id)= R(\phi)$ is the Reidemeister number of $\phi$ . If both  $\phi$ and $\psi$  are the identity maps then the $(id,id)$-conjugacy classes are the usual
conjugacy classes in the group $G$.
We shall write $\{x\}_{\phi,\psi}$ for the $(\phi,\psi)$-{\em conjugacy} or
{\em bitwisted conjugacy} class
 of the element $x\in G$.
\end{dfn}

If $G$ is a finite group, then the classical Burnside-Frobenius  theorem (see e.g.\cite{s1} )
says that the number of
classes of irreducible representations is equal to the number of conjugacy
classes of elements of $G$.  Let $\wh G$ be the {\em unitary dual} of $G$,
i.e. the set of equivalence classes of unitary irreducible
representations of $G$.

If $\phi: G\to G$ is an automorphism, it induces a map $\wh\phi:\wh G\to\wh G$,
$\wh\phi (\r)=\r\circ\phi$.
Therefore, by the Burnside-Frobenius  theorem, if $\phi$ and $\psi$ are the identity automorphisms
of any finite group $G$, then we have
 $R(\phi,\psi)=\#\Coin(\wh\phi,\wh\psi)$, where $\Coin(\wh\phi,\wh\psi)$ is
the set of coincidence points of $ \wh\phi$ and $\wh\psi$.

{This statement remains true for $\phi\ne\Id$, $\psi = Id$ and finite $G$ \cite{fh} . Indeed,
consider an automorphism $\phi$ of a finite group $G$. Then $R(\phi)$
is equal to the dimension of the space of twisted invariant functions on
this group. Hence, by Peter-Weyl theorem (which asserts the existence of
a two-side equivariant isomorphism
$C^*(G)\cong \bigoplus_{\r\in\wh G} \End(H_\r)$),  $R(\phi)$  is identified with the sum of dimensions
$d_\r$ of twisted invariant elements of $\End(H_\r)$, where $\r$ runs over
$\wh G$, and the space of representation $\r$ is denoted by $H_\r$.
By Schur lemma,
$d_\r=1$, if $\r$ is a fixed point of $\wh\phi$, and is zero otherwise.
Hence,
$R(\phi)$ coincides with the number of fixed points of $\wh\phi$.
}

The attempts to generalize Burnside-Frobenius   theorem to the case of non-identical
endomorphisms  and of non-finite group were inspired by the dynamical
questions and were the subject  of a series of papers
\cite{f1, fh, ft1, ftv, ft2}

\begin{rk}
If $\phi: G\to G$ is an epimorphism, it induces a map $\wh\phi:\wh G\to\wh G$,
$\wh\phi (\r)=\r\circ\phi$
(because a representation is irreducible if and only if the  scalar operators in the space of
representation are the only ones which commute with all operators of the
representation). This is also true for a general endomorphism of Abelian group $G$, but  is not the case for a general endomorphism $\phi$ of any  group.
\end{rk}

The paper consists of five Sections.
To make the presentation more detailed and  transparent we start  in Section 2 from a new approach for compact   groups and two automorphisms.
We develop this approach in Sections  3 and 4. 
After some preliminary and technical considerations we prove in Section 3
 the main
result  of the paper, namely

\begin{teo} {\sc Bitwisted Burnside-Frobenius  theorem}:
Let $G$ be an Abelian  group and $\phi$ and $\psi$ its endomorphisms.
 Then
$
R(\phi,\psi)=\#\Coin(\wh\phi,\wh\psi)
$
if one of these  numbers  is finite.
\end{teo}

  If $\psi$ is the identity map then bitwisted Burnside-Frobenius theorem
for Abelian groups
implies the twisted Burnside-Frobenius  theorem for Abelian groups from \cite{ft1}:
$R(\phi)=\#\Fix(\wh\phi)$
(the Reidemeister
number of an endomorphism $\phi$
is equal to the number of  fixed points of $\widehat\phi$
on the unitary dual, if one of these numbers is finite).

In some sense our theorem  is a reply to a remark of J.-P.~Serre
\cite{s1}, p.34  that for compact infinite groups
an analogue of Burnside-Frobenius  theorem is not interesting: $\infty=\infty$.
It turns out that for infinite discrete groups the situation differs
significantly, and even in non-twisted situations the number of
classes can be finite (for one of the first examples see another
book of J.-P.~Serre \cite{s2}). A number of examples
of groups and automorphisms with finite Reidemeister numbers was
 studied in \cite{f1,  fhw, ftv,  gw}.

In   Section \ref{Dehn} a bitwisted Dehn conjugacy problem
is formulated and  an affirmative
answer to this  problem
for almost polycyclic  groups  is obtained.

Finally, in Section \ref{sec:BaumSolit},  based on  \cite{fg1},   it is proved
that for any two  injective endomorphisms $\phi,\psi$ of a Baumslag-Solitar
 group $B(1,n)$ the  Reidemeister number $R(\phi,\psi)$ is infinite.
This is ``opposite'' case for the  bitwisted Burnside-Frobenius theorem.
It is interesting  open problem to describe all  groups with  property
that for any two injective endomorphisms $\phi,\psi$  
  the  Reidemeister number $R(\phi,\psi)$ is infinite.

\medskip
The interest in bitwisted and twisted  conjugacy relations has its origins
in the Nielsen-Reidemeister coincidence and fixed point  theory (see, e.g.
\cite{f1, w}, in Selberg theory (see, eg. \cite{ac},
and  Algebraic Geometry (see, e.g. \cite{g}).

The Lefschetz coincidence theorem generalizes the celebrated fixed point theorem to coincidences of two maps between closed connected orientable manifolds
of the same dimension. H. Schirmer \cite{sch} developed a Nielsen-Reidemeister  coincidence  theory and proved a Wecken type theorem for coincidences.The analogous
Lefschetz, Reidemeister and Nielsen coincidence numbers play a key role
in the Nielsen-Reidemeister  theory.

\medskip
\noindent
{\bf Acknowledgement.}
 I  would like to thank the MPI for its kind support and
hospitality while the most part of this work has been completed.

The author is  grateful to
D.~Gon{\c{c}}alves,
R.~Hill,
E.~ Troitsky,
L.~Vainerman and
A.~Vershik
for helpful discussions.

The author is also indebted to the referee for a careful reading and valuable
suggestions.

\section{Compact groups} \label{subsec:compcase}
Let $G$ be a compact  group and $\phi$ and $\psi$ its automorphisms.  Hence $\wh G$ is a discrete space. Then
 $C^*(G)=\oplus M_i$, where $M_i$ are the matrix algebras of
irreducible representations.
The infinite sum is in the following sense:
$$
C^*(G)=\{ f_i\}, i\in \{ 1,2,3,...\}=\hat G, f_i\in M_i,
\| f_i \| \to 0 (i\to \infty).
$$
When $G$ is finite and $\wh G$ is finite this is exactly Peter-Weyl theorem.

A characteristic function of a bitwisted class
is a functional on $C^*(G)$. For a  finite group it is evident,   for a
general compact group  it is necessary to verify only the measurability of
the bitwisted class
with the respect to Haar measure, i.e. that bitwisted class is Borel.
For a compact $G$,
the bitwisted conjugacy classes being orbits of bitwisted action are compact and
hence closed.
Then its complement is open, hence Borel, and the class is Borel  too.

Under the identification it passes to a sequence
$\{ \f_i \} $, where $\f_i$ is a functional on $M_i$ (the properties
of convergence can be formulated, but they play no role at the moment).
The conditions of invariance are the following: for each $\r_i \in\wh G$
one has $g[\f_i]=\f_i$, i.e. for any $a\in M_i$ and any $g\in G$ one has
 $\f_i(\r_i(\psi(g)) a \r_i(\phi(g^{-1})))=\f_i(a)$.

Let us recall  the following well-known fact.

\begin{lem}\label{lem:funcnamat}
Each functional on matrix algebra has form
$a\mapsto \Tr(ab)$ for a fixed matrix $b$.
\end{lem}

\begin{proof}  One has $\dim (M(n,\C))'=\dim (M(n,\C))=n\times n$
and looking
at matrices as at operators in $V$, $\dim V=n$, with base $e_i$,
one can remark that
functionals $a\mapsto \<a e_i, e_j\>$, $i,j=1,\dots, n$,
are linearly  independent. Hence,
any functional takes form
$$
a\mapsto \sum_{i,j} b_{ij} \<a e_i, e_j\>=
\sum_{i,j} b_{ij} a^j _i = \Tr(ba),\qquad
b:=\|b_{ij}\|.
$$

\end{proof}

\begin{teo}\label{teo:abelian}
Let $G$ be a compact   group and $\phi$ and $\psi$ its automorphisms.
 Then

$
R(\phi,\psi)=\#\Coin(\wh\phi,\wh\psi)
$
if one of these  numbers  is finite.
\end{teo}

\begin{proof}
 Let us consider  invariant functionals on matrix algebras under
bitwisted action of group $G$ :
$$
\Tr(b\r_i(\psi(g)) a \r_i(\phi(g^{-1})))=\Tr(ba),\qquad \forall\,a,g,
$$
$$
\Tr((b-\r_i(\phi(g^{-1})) b \r_i(\psi(g)))a)=0,\qquad \forall\,a,g,
$$
hence,
$$
b-\r_i(\phi(g^{-1})) b \r_i(\psi(g))=0,\qquad \forall\,g.
$$
Since $\r_i$ is irreducible, the dimension of the space of such $b$ is
1 if $\r_i$ is a coincidence
point of $\wh\phi$ and  $\wh\psi$ and 0 in the opposite case. So, we are done.

\end{proof}

\begin{rk}
In fact we are only interested in finite group  case.
Indeed, for a compact $G$,
the bitwisted conjugacy classes being orbits of bitwisted action are
compact and hence closed.
If there is a finite number of them, then are open too.
Hence, the situation is more or less
reduced to a discrete group: quotient by the component of unity.
\end{rk}
\begin{rk}
Bitwisted Burnside-Frobenius theorem for compact groups  remains  true if $\phi$ is automorphism
and $\psi$ is any endomorphism. It is equivalent to the twisted Burnside-Frobenius theorem for  the one  endomorphism
$\phi^{-1}\circ \psi$ (see \cite{ft2}).
\end{rk}
\begin{ex}
The following counterexample to the  bitwisted Burnside-Frobenius theorem  for two {\it endomorphisms} of a nonabelian finite group was proposed to the author  by E. Troitsky \cite{t1}.
Let us consider two trivial endomorphisms $\phi$ and $\psi$ of nonabelian finite group $G$, i.e.  $\phi(g)=\psi(g)=e$ for every element $g\in G$.
Then for every $x\in G$ and every $\gamma \in G$ we have $\psi(\gamma)x\phi(\gamma^{-1})=x$, hence  $R(\phi,\psi)=\#G$. From another side
$\wh\phi (\r)=\r\circ\phi = I =\r\circ\psi=\wh\psi (\r)$ for every $\r \in \wh G$  hence, if $G$ is nonabelian, then
$\# Coin(\wh\phi,\wh\psi)= \# \wh G$ and  $R(\phi,\psi)=\#G\neq\# \wh G= \# Coin(\wh\phi,\wh\psi)$.
\end{ex}

\section{ Endomorphisms of Abelian groups}\label{subsec:abelcase}

Let $\phi,\psi$ be   endomorphisms  of a  Abelian group $G$.

\begin{lem}\label{lem:charabelclassred}
Let $G$ be abelian.
The bitwisted conjugacy class $H$ of $e$ is a subgroup.
The other ones  are cosets $gH$.
\end{lem}

\begin{proof}
The first statement follows from the equalities
$$
(\psi(h)\phi(h^{-1})\psi( g)\phi(g^{-1})=\psi(gh) \phi((gh)^{-1},
$$
$$
(\psi(h)\phi(h^{-1}))^{-1}=\phi(h)\psi( h^{-1})=\psi(h^{-1}) \phi(h).
$$
For the second statement suppose $a\sim b$, i.e. $b=\psi(h)a\phi(h^{-1})$. Then
$$
gb=g\psi(h)a\phi(h^{-1})=\psi(h)(ga)\phi(h^{-1}), \qquad gb\sim ga.
$$
\end{proof}
\begin{lem}
Suppose, $u_1,u_2\in G$, $\chi_H$ is the
characteristic function of $H$ as a set. Then
$$
\chi_H(u_1 u_2^{-1})=\left\{
\begin{array}{ll}
1,& \mbox{ if } u_1,u_2 \mbox{ are in one coset },\\
0, & \mbox{ otherwise }.
\end{array}
\right.
$$
\end{lem}

\begin{proof}
Suppose, $u_1\in g_1 H$, $u_2\in g_2 H$, hence, $u_1=g_1 h_1$, $u_2=g_2 h_2$.
Then
$$
u_1 u_2^{-1}=g_1 h_1 h_2^{-1} g_2^{-1} \in g_1 g_2^{-1} H.
$$
Thus, $\chi_H(u_1 u_2^{-1})=1$ if and only if
$g_1 g_2^{-1}\in  H$ and $u_1$ and $u_2$
are in the same class. Otherwise it is 0.

\end{proof}
The following Lemma is well known.

\begin{lem}
For any subgroup $H$ the function $\chi_H$ is of positive type.
\end{lem}

\begin{proof}
Let us take arbitrary elements $u_1,u_2,\dots, u_n$ of $G$.
Let us reenumerate them
in such a way that some first are in $g_1 H$, the next ones are  in $g_2 H$,
and so on, till
$g_m H$, where $g_j H$ are different cosets. By the previous Lemma the matrix
$\|p_{it}\|:=\|\chi_H(u_i u_t^{-1})\| $ is block-diagonal with
square blocks formed by
units. These blocks, and consequently  the whole  matrix are  positively
semi-defined.

\end{proof}
\begin{lem}
In the Abelian case characteristic functions of bitwisted conjugacy classes
belong to
the Fourier-Stieltjes algebra $B(G)=(C^*(G))^*$.
\end{lem}

\begin{proof}
By Lemma \ref{lem:charabelclassred}
in this case the characteristic functions of bitwisted conjugacy
classes are the shifts
of the characteristic function of the class $H$ of $e$.
Hence, by Corollary (2.19) of [5], these characteristic functions are in
$B(G)$.

\end{proof}
Let us remark that there exists a natural isomorphism (Fourier transform)
$$
u\mapsto\wh u,\qquad C^*(G)=C_r^*(G)\cong C(\wh G),\qquad
\wh g(\rho):=\rho(g),
$$
(this is a  number because irreducible representations of an Abelian group
are 1-dimen\-sional). In fact, it is better to look (for what follows) at the  algebra $C(\wh G)$ as an algebra of continuous sections  of a bundle of 1-dimensional matrix algebras
over $\wh G$.

Our characteristic functions, being in $B(G)=(C^*(G))^*$ in this case,
are mapped to the functionals on $C(\wh G)$  which, by the
Riesz(-Markov-Kakutani)  theorem,
are measures on $\wh G$ \cite{c}. Which of these measures  are invariant under
the induced (bitwisted) action of $G$ ? Let us remark,
that an invariant non-trivial
functional gives rise to at least one invariant space -- its kernel.

Let us remark, that convolution under the Fourier transform becomes
point-wise multiplication. More precisely, the bitwisted action, for example,
is defined as
$$
g[f](\rho)=\rho(\psi(g))f(\rho)\rho(\phi(g^{-1})),\qquad \rho\in\wh G,\quad
g\in G,\quad f\in C(\wh G).
$$

There are 2 possibilities  for the bitwisted action of $G$ on the representation
 algebra $A_\rho \cong \C$: 1) the linear span of the orbit of $1 \in A_\rho$
is equal to all $A_\rho$, 2) and  the opposite case (the action is trivial).

The second case means that the space of intertwining operators between
 $A_{\wh \psi \rho}$ and $A_{\wh \phi \rho}$ equals $\C$, and  $\rho$ is a coincidence  point of
 $\wh \phi$ and $\wh \psi$. In the first case this is the opposite
situation.

If we have a finite number of such coincidence
points, then the space of bitwisted invariant measures is just the space of
measures
concentrated in these points. Indeed, let us describe the action of $G$
on measures in
more detail.

\begin{lem}
For any  Borel set $E$ one has
$g[\mu](E)=\int_E g[1]\,d\mu$.
\end{lem}

\begin{proof}
The restriction of measure to any Borel set commutes with the action of  $G$,
since the last is point wise on $C(\wh G)$.
For any  Borel set $E$ one has
$$
g[\mu](E)= \int_E 1\, dg[\mu]= \int_E g[1] \,d\mu.
$$

Hence, if $\mu$ is bitwisted invariant, then for any  Borel set $E$ and any
$g\in G$ one has
$$
\int_E (1-g[1])\,d\mu =0.
$$
\end{proof}
\begin{lem}
Suppose, $f\in C(X)$, where $X$ is a compact Hausdorff space, and $\mu$
is a regular Borel measure on $X$, i.e. a functional on $C(X)$.
Suppose, for any
Borel set $E\subset X$ one has $\int_E f\, d\mu=0$.
Then $\mu (h)=0$ for any $h\in C(X)$
such that $f(x)=0$ implies $h(x)=0$. I.e. $\mu$ is concentrated off the
interior
of $\rm supp\,f$.
\end{lem}

\begin{proof}
Since the functions of the form $fh$ are dense in the space of functions   $h$'s,
it is sufficient to verify the statement for $fh$. Let us choose an
arbitrary $\e>0$
and a simple function $h'=\sum\limits_{i=1}^n a_i \chi_{E_i}$
such that $|\mu(fh')-\mu(fh)|<\e$.
Then
$$
\mu(fh')=\sum_{i=1}^n \int_{E_i} a_i f \,d\mu =\sum_{i=1}^n a_i
\int_{E_i}f \,d\mu =0.
$$
Since $\e$ is an arbitrary one, we are done.
\end{proof}

Applying this lemma to a bitwisted invariant measure $\mu$ and $f=1-g[1]$
we obtain that $\mu$ is concentrated at our finite number of coincidence
points of $\wh\phi$ and $\wh \psi$, because outside of them $f\ne 0$.

If we have an infinite number of coincidence  points, then
 the space is
infinite 
dimensional(we have an infinite number of measures
concentrated in finite number of points, each time different)
and Reidemeister number is infinite as well. So, we have  proved the
following theorem:

\begin{teo}\label{teo:abelian}
Let $G$ be an Abelian  group and $\phi$ and $\psi$ its endomorphisms.

 Then
$
R(\phi,\psi)=\#\Coin(\wh\phi,\wh\psi)
$
if one of these  numbers  is finite.
\end{teo}

\section{ Automorphisms of bitwisted conjugacy separable groups, weak bitwisted
Burnside-Frobenius theorem}

 If  the endomorphisms  $\phi$ and $\psi$ are the  automorphisms,  then  the 
bitwisted Burnside-Frobenius problem is equivalent to the twisted Burnside-Frobenius  problem for the one  automorphism
$\phi^{-1}\circ \psi$. So we can apply all the  results  from \cite{ft2, ft3,  t2} for the one automorphism. 
In the present paper we define  the following property  for a countable
discrete group $G$
\begin{dfn}
A group $G$ is called \emph{$(\phi,\psi)$-conjugacy separable} with respect to
an endomorphisms $\phi,\psi:G\to G$ or \emph{bitwisted conjugacy separable}, if any pair $g$, $h$ of
   non-$(\phi,\psi)$-conjugate elements of $G$ are non-$(\ov\phi,\ov\psi)$-conjugate in some finite
   quotient of $G$ respecting $\phi,\psi$.
\end{dfn}
 If $\psi=id$ and $\phi$ is  automorphism
then this definition gives us  the definition of $\phi$-conjugacy separable discrete  group
 \cite{ft2}. If both  $\phi=id$ and $\psi=id$ then this definition gives us well known notion of conjugacy separable group.

Now  we formulate  the main 
results  for two automorphisms, namely

\begin{enumerate}
    \item { $(\phi,\psi)$- conjugacy separable groups and  some extensions:}
Suppose, there is an extension $H\to G\to G/H$, respecting automorphisms  $\phi$ and $\psi$ where the normal invariant subgroup $H$ is a
{\rm $(\phi',\psi')$- conjugacy separable} group  ;
 Reidemeister number  $R(\phi,\psi)$ is finite;
 $G/H$ is finitely generated {\rm FC}-group (i.e.
a group with finite conjugacy classes).
Then $G$ is an $(\phi,\psi)$- conjugacy separable group ).
    \item { Classes of $(\phi,\psi)$- conjugacy  separable  groups:} almost polycyclic groups and
    finitely generated groups of polynomial growth are $(\phi,\psi)$- conjugacy  separable groups if   $\phi$ and $\psi$ are automorphisms.
   
    \item { Bitwisted Burnside-Frobenius  theorem:}
Let $G$ be an {\rm $(\phi,\psi)$- conjugacy separable} group and $\phi$ and $\psi$ its automorphisms.
Denote by $\wh G_f$ the subset of the unitary dual $\wh G$ related to
finite-di\-men\-si\-on\-al representations and
by $C_f(\phi,\psi)$
the number of coincidence  points of $\wh\phi_f$ and $\wh\psi_f$ on $\wh G_f$. Then
$
R(\phi,\psi)=C_f(\phi,\psi)
$
if one of these  numbers  is finite .

\item { The affirmative answer to the bitwisted Dehn conjugacy problem for
polycyclic-by-finite groups} ( see subsection \ref{Dehn}).
\end{enumerate}

  If $\psi$ is the identity map then bitwisted Burnside-Frobenius theorem
implies the twisted Burnside-Frobenius  theorem from \cite{ft2}:
$R(\phi)=\#\Fix(\wh\phi_f)$
(the Reidemeister
number of an automorphism $\phi$
is equal to the number of finite-dimensional fixed points of $\widehat\phi$
on the unitary dual, if one of these numbers is finite).

 In \cite{ft1} one obtained  from the
twisted Burnside-Frobenius  theorem the congruences for Reidemeister numbers which
are similar to the remarkable Dold congruences for the Lefschetz numbers
and which,  together with the twisted Burnside-Frobenius  theorem itself,
is very important for the realization problem of Reidemeister numbers in
topological dynamics and the study of the Reidemeister
zeta-function.

We would like to emphasize the following important
remarks.

\begin{enumerate}
\item {In the original formulation by Fel'shyn and Hill \cite{fh}
the conjecture about twisted Burnside-Frobenius theorem asserts an equality
of $R(\phi)$ and the number of fixed points of $\wh\phi$ on $\wh G$. This
conjecture was proved in \cite{fh,ft1} for f.g. type I groups.}

\item As it follows from a key example, which we have studied  in
\cite{ftv}, a twisted conjugacy separable group can have infinite
dimensional "supplementary" coincidence  representations. More precisely
we discuss in that paper the case of a semi-direct product of the action
of $\Z$ on $\Z\oplus \Z$ by a hyperbolic automorphism with finite
Reidemeister number (four to be precise) and the number of fixed points
of $\wh\phi$ on $\wh G$ equal or greater than five, while the number
of fixed points on $\wh G_f$ is four.

{This gives a counterexample to the conjecture in its original
formulation and leads to the formulation using  only fixed(or  coincidence)
points in
$\wh G_f$. }

\item
The origin of the
phenomenon  of an extra coincidence  point
lies in bad separation properties of $\wh G$ for general
discrete groups. A more deep study leads  to the following general theorem.

 \begin{teo} {\sc Weak bitwisted Burnside-Frobenius  theorem}:
 Let $G$ be a countable discrete  group and $\phi$ and $\psi$ its automorphisms. The number $R_*(\phi,\psi)$ of $(\phi,\psi)$- conjugacy classes related to
bitwisted invariant functions on $G$ from the Fourier-Stieltjes
algebra $B(G)$ is  equal to the number $ C_*(\phi,\psi)$  of generalized
coincidence  points of $\wh\phi$ and  $\wh\psi$ on the Glimm spectrum of $G$, i.~e. on the complete regularization of $\wh G$, if one of the numbers  $R_*(\phi,\psi)$ and $ C_*(\phi,\psi)$  is finite.
 \end{teo}

\item
Nevertheless for extreme  groups( like Osin group \cite{ft1}, which is infinite finitely generated simple  group  with two conjugacy classes),
even the modified conjecture is not true.
Keeping in mind that for non-elementary Gromov hyperbolic 
 groups, mapping class groups, saturated weakly branch groups and any  two  automorphisms
$\phi,\psi$ the Reidemeister number  $R(\phi,\psi)$ is  infinite
\cite{f2, ll, fg2, flt}  as well as for any two injective  endomorphisms of
Baumslag-Solitar groups 
cf. Section \ref{sec:BaumSolit} and \cite{fg1}
while in the "opposite" case the bitwisted Burnside-Frobenius theorem is proved, we
can hope that our results   can lead
to a complete resolution of the problem, if the extreme groups
will be handled.

\end{enumerate}

Let us remark, that in some cases the weak bitwisted Burnside-Frobenius  theorem
easily implies the bitwisted Burnside-Frobenius  theorem (in particular, in the form with
finite-dimensional representations). For example, we will show
directly that $R(\phi,\psi)=R_*(\phi,\psi)$ in Abelian case. On the other
hand, the unitary dual coincides with the Glimm spectrum.

\subsection{Bitwisted Dehn conjugacy problem}\label{Dehn}

The subject is closely related to some decision problem. Recall that
M.~Dehn in 1911 \cite{d}
has formulated in particular

\smallskip\noindent
\textbf{Conjugacy problem:} Does there exists an algorithm to determine
whether an arbitrary pair of group words $U$, $V$ in the generators of $G$ define
conjugate elements of $G$?

\smallskip
The following question was posed by G.~Makanin \cite{m2} [Question 10.26(a)]:

\smallskip\noindent
\textbf{Question:} Does there exists an algorithm to determine
whether for an arbitrary pair of group words $U$ and $V$
of a free group $G$ and an arbitrary automorphism $\phi$ of $G$
the equation $\phi(X)U=VX$ solvable in $G$?

\smallskip
In  \cite{b2} an  affirmative answer  to the Makanin's question
is obtained.

\smallskip
In  \cite{b1}  the following problem, which generalizes the two
above problems, was posed:

\smallskip\noindent
\textbf{Twisted conjugacy problem:} Does there exists an algorithm to determine
whether for an arbitrary pair of group words $U$ and $V$ in the generators of $G$
the equality $\phi(X)U=VX$ holds for some $W\in G$ and $\phi\in H$, where $H$ is
a fixed subset of $\Aut(G)$?

In \cite{ft4}  we proved that the twisted conjugacy problem has the affirmative
answer for $G$ being polycyclic-by-finite group and $H$ be equal to a unique automorphism $\phi $.

In this article we suggest the following problem, which generalises the three above problems:

\smallskip\noindent
\textbf{Bitwisted conjugacy problem:} Does there exists an algorithm to determine
whether for an arbitrary pair of group words $U$ and $V$ in the generators of $G$
the equality $\psi(X)U=V\phi(X)$ holds for some $W\in G$ and endomorphisms $\phi,\psi\in H$, where $H$ is
a fixed subset of $\End(G)$?

\begin{teo}\label{teo:twistconjprobforpolbfin}
The bitwisted conjugacy problem has the affirmative answer for $G$ being
almost polycyclic group and $H$ be equal to a pair of  automorphisms $\phi,\psi$.
\end{teo}

\begin{proof}
It follows immediately from  $(\phi,\psi)$- conjugacy  separability of  almost polycyclic groups  by the same
argument as in the paper of
Mal'cev \cite{m1} (see also \cite{m3}, where the property
of conjugacy separability was first formulated)
for the (non-twisted) conjugacy problem.

\end{proof}

\section{Baumslag-Solitar groups $B(1,n)$}\label{sec:BaumSolit}

In this section based on (a part of) [8]  it is proved
that for any two injective endomorphisms $\phi,\psi$ of a Baumslag-Solitar group $B(1,n)$ the  Reidemeister number $R(\phi,\psi)$ is infinite.
 From  \cite{fg1} the similar result for two  injective endomorphisms  of $B(m,n)$  where $|n|\ne |m|$ and $|nm|\ne 0,1$  follows  as well.
For this   class of groups  there
 is no  bitwisted Burnside-Frobenius theorem
( the Reidemeister number $R(\phi,\psi)$  is  infinite).
It is  open problem to describe all  groups with such property.

The presentation in this section  is selfcontained.

Let $B(1,n)=\langle a,b: a^{-1}ba=b^n, n>1 \rangle $ be
the Baumslag-Solitar groups. These groups 
are finitely-generated solvable groups (in particular amenable) which are not
virtually nilpotent. These  groups have
exponential growth   and  they are not Gromov hyperbolic.
Furthermore, these groups are  torsion  free and metabelian (an extension of
an Abelian group by an Abelian). More precisely one has

\begin{prop}\label{prop:felgon3.0}
$
B(1,n)\cong {\Z}[1/n]\rtimes_{\theta} {\Z},
$
where the
action of ${\Z}$ on ${\Z}[1/n]$ is given by $\theta(1)(x)=x/n$.
\end{prop}

\begin{proof}
The map defined by   $\iota(a)=(0,1)$ and
$\iota (b)=(1,0)$ extends to a  unique homomorphism
$\iota:B(1,n) \to \Z[1/n]\rtimes \Z,$
because
$$
\iota(a^{-1}) * \iota(b) * \iota(a)=
(0,-1) * (1,0) * (0,1)=(0,-1) * (1,1)=(n,0)=\iota(b^n).
$$
One has
\begin{equation}\label{ur:felgondob}
  \iota(a^r b^s a^{-r})=(0, r) * (s, -r)= \left(\frac
  s{n^r},0\right).
\end{equation}

The map $\iota$ is clearly surjective. Let us show
that this homomorphism is injective.
Let us remark that the group relation implies
$a^{-1}b^{-1}a=b^{-n}$. Hence for any $s$ one has
$a^{-1} b^s a=b^{ns}$.
Thus we can move all $a^{-1}$ to the right (until they
annihilate with some $a$ or take the extreme right place)
and all $a$ to the left, with an appropriate changing of powers
of $b$. Hence
any word in $B(1,n)$
is equivalent to a word of the form
$w=a^{r_1}b^{s}a^{r_2}$, where $r_1 \ge 0$, $r_2 \le 0$.
Then
$\iota(w)=(m,r_1 + r_2)$ for some $m\in \Z[1/n]$.
Hence, if $r_1+r_2\ne 0$, then $\iota(w)\ne e$. Let now
$r_1+r_2=0$, then by (\ref{ur:felgondob})
if $\iota(w)=e$, then $s=0$ and $w=e$.

Consider the homomorphism $|\ \ |_a: B(1,n) \longrightarrow \Z $ which
associates to each word $w\in B(1,n)$ the sum of the exponents of $a$ in
this word. Since this sum for the relation is zero,
this is a well defined map, which
is evidently surjective.
\end{proof}
\begin{prop}\label{prop:felgon3.1}\cite{fg1}
We have a short exact sequence
$$
0 \longrightarrow K\longrightarrow B(1,n)
\stackrel{|\ \ |_a}\longrightarrow \Z \longrightarrow 1,
$$
where $K$ is the kernel of $ |\ \ |_a$.
Moreover, $B(1,n)$ equals a semidirect product $K \rtimes \Z$.
\end{prop}

\begin{proof}
The first statement follows from surjectivity of $|\cdot|_a$.
Since $\Z$ is free, this sequence splits.
\end{proof}

\begin{prop}\label{prop:felgon3.2}\cite{fg1}
The kernel $K$ coincide with the
normalizer  $ N\langle b\rangle $ of the subgroup $\langle b \rangle$
generated by $b$ in $B(1,n)$:
\begin{equation}\label{ur:felgon3.2}
  0 \longrightarrow N\langle b\rangle \longrightarrow B(1,n)
\stackrel{|\ \ |_a}\longrightarrow \Z \longrightarrow 1.
\end{equation}
\end{prop}

\begin{proof}
We have $ N\langle b\rangle\subset K $. The quotient $B(1,n)/N\langle
b\rangle $ has the following presentation: $\ov a^{-1}\ov b \ov a=
\ov b^n$, $\ov b=1$. Therefore this group is isomorphic to ${\Z}$
 under the identification $[a]\leftrightarrow 1_\Z$.
Hence the
natural projection coincides with the map $|\ \ |_a$
and we have the following commutative diagram
$$
\xymatrix{
0 \ar[r] &  N\langle b\rangle \ar[r]\ar[d]^\cap & B(1,n) \ar[r]\ar[d]^= &
B(1,n)/(N\langle b\rangle)\ar[r]\ar[d]^\cong & 1  \\
0 \ar[r]& K \ar[r] & B(1,n) \ar[r] & {\Z} \ar[r] & 1.
}
$$
The five-lemma completes the proof.

\end{proof}
\begin{prop}\label{prop:felgon3.4}\cite{fg1}
Any endmorphism $\phi: B(1,n) \to B(1,n)$ is a endomorphism of
the short exact sequence {\rm (\ref{ur:felgon3.2})}.
\end{prop}

\begin{proof}
 Let $\bar \phi$ be the  endomorphism induced by $\phi$ on
the abelianization $B(1,n)_{ab}$ of $B(1,n)$. The group
$B(1,n)_{ab}$ is isomorphic to $Z_{n-1} +Z$. The torsion elements of
$B(1,n)_{ab}$ form a subgroup isomorphic to $Z_{n-1}$ which is invariant
under any endomorphism. The preimage of this subgroup under the projection
$B(1,n) \to B(1,n)_{ab}$ is  exactly the subgroup
$N(b)$, i.e. the elements represented by words where the sum of the powers
of $a$ is zero. So it follows that $N(b)$ is mapped into $N(b)$ and the
result follows.
\end{proof}

\begin{teo}\label{teo:felgon3.4}\cite{fg1}
For any two injective
endomorphisms $\phi$ and  $\psi$  of $B(1,n)$ the  Reidemeister number $R(\phi,\psi)$  is
infinite.
\end{teo}

\begin{proof}
By Proposition \ref{prop:felgon3.4} they  are  
endomorphisms of the short exact sequence. The induced maps  $\ov\phi$ and
$\ov\psi$
on the quotient are
 injective endomorphisms of $\Z$. If $\ov \phi=\Id_\Z$,
and $\ov \psi=\Id_\Z$,
then by \cite{fg1} the number of bitwisted 
classes is  infinite. Hence,at least one of  $\ov \phi$
or $\ov \psi$  is multiplication by $k\ne 0,1$.
But this is impossible. Indeed, when we apply
$\phi$ to the relation $a^{-1}ba=b^n$, under the identification of
Proposition \ref{prop:felgon3.0}
$\iota:B(1,n)\cong \Z[1/n]\rtimes \Z$  we have,
because $\phi(b)\in N\<b\>$ and hence $\iota(\phi(b))=(d,0)$ for
some $d\in\Z[1/n]$,
$$
(n d,0)=\iota(\phi (b^n))=
\iota(\phi(a^{-1}ba))=\iota(a^{-k} \phi(b) a^{k})=
(d\cdot n^k,0).
$$
This implies that either  $n^{1-k}=1$ or $\phi(b)=0$.
\end{proof}
\def\cprime{$'$} \def\dbar{\leavevmode\hbox to 0pt{\hskip.2ex \accent"16\hss}d}
  \def\polhk#1{\setbox0=\hbox{#1}{\ooalign{\hidewidth
  \lower1.5ex\hbox{`}\hidewidth\crcr\unhbox0}}}
\providecommand{\bysame}{\leavevmode\hbox to3em{\hrulefill}\thinspace}
\providecommand{\MR}{\relax\ifhmode\unskip\space\fi MR }
\providecommand{\MRhref}[2]{%
  \href{http://www.ams.org/mathscinet-getitem?mr=#1}{#2}
}
\providecommand{\href}[2]{#2}


\begin{thebibliography}{10}

\bibitem{ac}
 J. Arthur and L. Clozel,
 {\it Simple algebras, base change, and the advanced
  theory of the trace formula},
 Princeton University Press, Princeton, NJ,1989.

\bibitem{b1} 
Valerij Bardakov, Leonid Bokut, and Andrei Vesnin,
 {\it Twisted conjugacy in
  free groups and Makanin's question}, Southeast Asian Bull. Math.
  \textbf{29} (2005), no.~2, 209--226, (E-print arxiv:math.GR/0401349).

\bibitem{b2} 
O. Bogopolski, A. Martino, O. Maslakova, E. Ventura,
{\it Free-by-cyclic groups have solvable conjugacy problem}, Bulletin of the London
Mathematical Society 38(5)(2006), 787-794.  
\bibitem{c}
 A. Connes, {\it Noncommutative Geometry},
Academic Press, 1994.

\bibitem{d}
 M.~Dehn, 
{\it \"Uber unendliche diskontinuierliche Gruppen}, Math. Ann.
  \textbf{71} (1911), no.~1, 116--144. 

\bibitem{e}
 P.~Eymard, 
{\it L'alg\`ebre de Fourier d'un groupe localement compact},
  Bull. Soc. math. France \textbf{92} (1964), 181--236.

\bibitem{f1}
 A.~Fel'shtyn,
 {\it Dynamical zeta functions, Nielsen theory and
  Reidemeister torsion}, Mem. Amer. Math. Soc. \textbf{147} (2000), no.~699,
  xii+146.
\bibitem{f2} 
A.~L. Fel'shtyn, 
{\it The Reidemeister number of any automorphism of a
  Gromov hyperbolic group is infinite}, Zap. Nauchn. Sem. S.-Peterburg.
  Otdel. Mat. Inst. Steklov. (POMI) \textbf{279} (2001), no.~6 (Geom. i
  Topol.), 229--240, 250.

\bibitem{fg1}
 A.~Fel'shtyn and D.~Gon{\c{c}}alves,
{\it Reidemeister numbers of any automorphism of  Baumslag-Solitar group}, is infinite. Geometry and Dynamics of Groups and Spaces, Progress in Mathematics, v. 265(2008), 399 - 414.

\bibitem{fg2}
 Alexander Fel'shtyn and Daciberg L. Gon\c{c}alves,
     {\it  Twisted conjugacy classes in Symplectic groups, Mapping class groups and Braid groups(including an Appendix written  with Francois Dahmani)}.E-print arxiv:
math.GR/0708.2628.

\bibitem{fh} 
A.~Fel'shtyn and R.~Hill,
 {\it The Reidemeister zeta function with
  applications to Nielsen theory and a connection with Reidemeister
  torsion}, $K$-Theory \textbf{8} (1994), no.~4, 367--393. 

\bibitem{fhw}
 Alexander Fel'shtyn, Richard Hill, and Peter Wong,
 {\it Reidemeister
  numbers of equivariant maps}, Topology Appl. \textbf{67} (1995), no.~2,
  119--131.

\bibitem{flt}
 A.~Fel'shtyn,Y. Leonov, E.~Troitsky,
 {\it Reidemeister numbers of saturated weakly branch groups}.
Geometria Dedicata, v.134(2008), 61-73.

\bibitem{ft1}
 A.~Fel'shtyn and E.~Troitsky,
 {\it A twisted Burnside theorem for countable
  groups and Reidemeister numbers}, Proc. Workshop Noncommutative Geometry
  and Number Theory (Bonn, 2003) (K.~Consani, M.~Marcolli
  eds.), Vieweg, Braunschweig, 2006, pp.141 - 154.

\bibitem{ftv} 
A.~Fel'shtyn, E.~Troitsky, and A.~Vershik,
 {\it Twisted Burnside theorem for
  type {II}${}_1$ groups: an example}, Mathematical Research Letters 13(2006), no.5, 719-728.

\bibitem{ft2}
  A.~Fel'shtyn and E.~Troitsky, 
{\it  Twisted {B}urnside-Frobenius  theory for discrete groups} 
J. Reine Angew. Math.(Crelle's Journal) 613(2007), 193-210.
 

\bibitem{ft3} 
A.~Fel'shtyn and E.~Troitsky, 
{\it Geometry of Reidemeister classes
 and twisted Burnside theorem} Journal of K-theory,  vol 2, issue 3(2008), 405-445.

\bibitem{ft4}
 A.~Fel'shtyn and E.~Troitsky,
 { \it Twisted conjugacy separable groups},

E-print arXiv:math.GR/0606764, 2006.

\bibitem{gw} D.~Gon{\c{c}}alves and P.~Wong,
 {\it Twisted conjugacy classes in exponential
  growth groups}, Bull. London Math. Soc. \textbf{35} (2003), no.~2, 261--268.
  
\bibitem{g} A.~Grothendieck,
 {\it Formules de {N}ielsen-{W}ecken et de {L}efschetz en
  g\'eom\'etrie alg\'ebrique}, S\'eminaire de G\'eom\'etrie Alg\'ebrique du
  {B}ois-{M}arie 1965-66. {SGA} 5, Lecture Notes in Math., vol. 569,
  Springer-Verlag, Berlin, 1977, pp.~407--441.

\bibitem{ll}
 G.~Levitt and M.~Lustig,
{\it {Most automorphisms of a hyperbolic group have
  very simple dynamics.}}, Ann. Scient. \'Ec. Norm. Sup. \textbf{33} (2000),
  507--517.

\bibitem{m1}
 A.I. Malcev,
 {\it On homomorphisms onto finite groups}, Uchen. Zapiski
  Ivanovsk. ped. instituta \textbf{18} (1958), no.~5, 49--60, (=Selected
  Papers, Vol.~1, 1976, 450--461).

\bibitem{m2} 
V.~D. Mazurov and E.~I. Khukhro (eds.), 
{\it The {K}ourovka notebook},
  augmented ed., Ros. Akademiya Nauk Sibirskoe Otdelenie, Institut
  Matematiki im. S. L. Soboleva, Novosibirsk, 2002, Unsolved problems in group
  theory.

\bibitem{m3}
 A.~W. Mostowski,
 {\it On the decidability of some problems in special classes
  of groups}, Fund. Math. \textbf{59} (1966), 123--135. 

\bibitem{sch} H.~Schirmer,
{\it Mindestzahlen von Koinzidenzpunkten.}
J. reine angew. Math. 194(1955), 21-39.

\bibitem{s1} Jean-Pierre Serre,
 {\it Linear representations of finite groups},
  Springer-Verlag, New York, 1977, Translated from the second French edition by
  Leonard L. Scott, Graduate Texts in Mathematics, Vol. 42. 

\bibitem{s2}
  Jean-Pierre Serre,
{\it Trees}, Springer Monographs in Mathematics, Springer-Verlag,
  Berlin, 2003, Translated from the French original by John Stillwell,
  Corrected 2nd printing of the 1980 English translation.
 
\bibitem{t1} E.~Troitsky,
{\it Private communication}.

\bibitem{t2} E.~Troitsky,
 {\it Noncommutative Riesz  theorem and weak
   Burnside  type theorem on twisted conjugacy}, Funct. Anal. Pril. \textbf{40} (2006), no.~2,
  44--54,

\bibitem{w} P.~ Wong,
{\it Reidemeister number, Hirsch rank, coincidences on polycyclic groups and
solvmanifolds.}
J. reine angew. Math. 524(2000), 185-204.

\end{thebibliography}
\end{document}